\documentclass[12pt,leqno]{article}

\usepackage{amssymb}
\usepackage{amsmath}
\usepackage{color}
\usepackage{url}
\newcommand{\comment}[1]{}
\newcommand{\de}{\mbox{\ :=\ }}

\newcommand{\Mm}{\mbox{$\mathfrak M$}}

\newcommand{\Nn}{\mbox{$\mathfrak N$}}
\newcommand{\Pp}{\mbox{$\mathfrak P$}}

\newcommand{\Gg}{\mbox{$\mathfrak G$}}
\newcommand{\Hh}{\mbox{$\mathfrak H$}}

\newcommand{\Aa}{\mbox{$\mathfrak A$}}
\newcommand{\Bb}{\mbox{$\mathfrak B$}}
\newcommand{\Aut}{\mbox{\sf Aut}}
\newcommand{\AutSpec}{\mbox{\sf AutSpec}}

\newcommand{\Set}{\mbox{\sf Set}}

\newcommand{\Th}{\mbox{\sf Th}}
\newcommand{\T}{\mbox{\sf T}}

\newcommand{\K}{\mbox{\sf K}}
\newcommand{\Mod}{\mbox{\sf Mod}}
\newcommand{\Cat}{\mbox{$\mathcal Mod$}}

\newcommand{\R}{\mbox{\sf R}}
\newcommand{\RR}{\mathcal{R}}
\newcommand{\Es}{\mbox{\sf S}}
\newcommand{\N}{\mathbb N}
\newcommand{\Z}{\mathbb Z}
\newcommand{\G}{\mbox{\sf G}}
\newcommand{\cG}{\mathcal G}

\newcommand{\zero}{\mbox{\sf 0}}
\newcommand{\suc}{\mbox{\sf suc}}

\newcommand{\qed}{\hfill\mbox{\bf QED}\bigskip}

\newtheorem{thm}{Theorem}
\newtheorem{lem}{Lemma}

\newtheorem{cor}{Corollary}

\comment{
\title[definitional equivalence via automorphism groups]{Testing definitional equivalence of theories via automorphism groups}
\author[H. Andr\'eka, J. Madar\'asz, I. N\'emeti and G. Sz\'ekely]
{HAJNAL ANDR\'EKA, JUDIT MADAR\'ASZ, ISTV\'AN N\'EMETI\\
\affil{Alfr\'ed R\'enyi Institute of Mathematics\\
	Budapest, Re\'altanoda st.\ 13-15, H-1053 Hungary
	\and GERGELY SZ\'EKELY%
	\thanks{This research is supported by the Hungarian National Research, Development and Innovation Office (NKFIH), grant no. FK-134732.}\\
\affil{Alfr\'ed R\'enyi Institute of Mathematics\\
	Budapest, Re\'altanoda st.\ 13-15, H-1053 Hungary
	\\[2.5pt]
	{\rm and}\\[-5pt]
	University of Public Service\\
	Budapest, 2 Ludovika tér, H-1083 Hungary
}}}}

\begin{document}
	
	\title{Testing definitional equivalence of theories via automorphism groups}
	\author{Andr\'eka, H., Madar\'asz, J., N\'emeti, I.\ and Sz\'ekely, G.%
		\thanks{This research is supported by the Hungarian National Research, Development and Innovation Office (NKFIH), grant no. FK-134732.}}
	\maketitle

\begin{abstract} 
Two first-order logic theories are definitionally equivalent if and only if there is a bijection between their model classes that preserves isomorphisms and ultraproducts (Theorem 2). This is a variant of a prior theorem of van Benthem and Pearce. In Example 2, uncountably many pairs of definitionally inequivalent theories are given such that their model categories are concretely isomorphic via bijections that preserve ultraproducts in the model categories up to isomorphism. Based on these results, we settle several conjectures of Barrett, Glymour and Halvorson.

\end{abstract}

\section{Introduction}

{\bf Classical definitional equivalence.} The subject of the present paper is the notion of (classical)
definitional equivalence of first-order logic theories. There are
various definitions of this notion scattered in the literature. Most
of these define the notion for theories with disjoint languages
only. We use the version defined in Lefever and Sz\'ekely \cite[Definition 11]{LSz} which
does not require the languages to be disjoint. According to this
definition, definitional equivalence of theories is the symmetric
and transitive closure of the relation ``definitional extension". 
This notion of definitional equivalence is shown to be the same as
the more prevailing ones for disjoint languages. For example, it coincides
with inter-translatability (\cite[Theorem 8]{LSz}) and ``having a
joint definitional extension" (\cite[Theorem 4]{LSz}). We believe
that making the vocabularies of theories disjoint is a
superfluous administrative task. Besides, making vocabularies disjoint
masks important intuitive features in many cases. This would be the case
in the present paper, too, e.g., in Example 2 and Theorem \ref{postest-cor}.

Definitional equivalence is also defined by means of a bijection
between two model classes in Henkin, Monk, and Tarski \cite[p.56]{HMT71}.
According to this definition, two theories are definitionally equivalent 
 when there is a bijection between their model classes such that connected 
 models are definitionally equivalent via the same definitions. This property is
called ``model mergeability" in \cite[Definition 13]{LSz} and is
proved to coincide with definitional equivalence as used in this
paper (\cite[Theorem 7]{LSz}). One of the advantages of
model mergeability is that it is kind of language-free in so far
that it is insensitive to whether the signatures of the two theories
overlap or not. Model mergeability is a mix of semantic and
syntactic features. 

A purely semantic characterization of definitional equivalence is given in de Bouv\`ere~\cite{Bou}, as follows. Two theories on disjoint languages are definitionally equivalent if and only if there is a third theory on the union of their languages such that both reduct-formation functions, from the model class of the third theory to the model classes of the two theories respectively, are bijections.
For variants of this characterization, see Barrett \cite[Corollary 2]{Bsym} and Lutz \cite[Claim 4]{Lutz}.  This semantic characterization is in terms of the concrete reduct-formation functions between model classes. Theorem~\ref{main-thm} in the present paper is a similar characterization for definitional equivalence: two theories are definitionally equivalent if and only if there is a bijection between their model classes that preserves universes, isomorphisms and ultraproducts. This is a purely semantic characterization of definitional equivalence similar to the one in \cite{HMT71} and different from the one in \cite{Bou}. The difference is that no third theory is used and arbitrary function is used in place of the concrete reduct-formation one. The idea of using functions that preserve isomorphisms and ultraproducts already occurs in van Benthem and Pearce \cite{vBP} where relative interpretability between first-order theories is characterized in place of definitional equivalence. For more on this, see Remark 5.
\bigskip
\goodbreak

{\bf Philosophy of science.} Definability theory is used quite extensively in recent philosophy of science papers, see for example \cite{BH22, Gly, Halvorson, HudetzDCE, Weatherall}. 
In philosophy of science, just as in mathematical logic, several notions of equivalence are used for comparing theories.  One is many-sorted definitional equivalence (\cite{AN, harnik, Mad02}) which is also called many-dimensional definitional equivalence (\cite{Hodges, Visser}) or Morita-equivalence (\cite{BH16, Halvorson}). Many-sorted definitional equivalence allows one to re-define the universes of models in a theory, therefore it is rather important. To distinguish definitional equivalence from many-sorted one, we sometimes call it classical definitional equivalence. Another version of equivalence of theories is bi-interpretability (see \cite{Hodges, Visser}). Categorical equivalence of theories (\cite{BH16, Weatherall}) is perhaps the weakest among the equivalences used for comparing theories. 

It is shown in Barrett and Halvorson \cite{BH16} that classical definitional equivalence, many-sorted definitional equivalence and categorical equivalence of theories are strictly weak\-er in this order.%
\footnote{It is not clear to us how bi-interpretability fits into this sequence.}
 Example 2 in this paper contains  pairs of theories on finite signatures that are categorically equivalent but not many-sorted definitionally equivalent (nor bi-interpretable). With this, we answer Barrett and Halvorson's questions \cite[Question 6.1]{BH16} and \cite[Question 1, p.77]{BDis17} concerning the importance of infinite signature in their counterexample.  In this context, it is natural to ask how much weaker categorical equivalence is than many-sorted definitional equivalence. Theorem \ref{main-thm} and especially its corollaries Theorem \ref{ultracat-thm} and Corollary \ref{dist-cor} in the present paper provide a property $\Pp$ of functors such that a functor establishing the categorical equivalence satisfies $\Pp$ if and only if the theories are classically definitionally equivalent. This property is that the functor is concrete and preserves ultraproducts. This is an answer to Barrett \cite[the question below Corollary 2]{Bsym}, \cite[Question 2]{BDis17} and Weatherall \cite[Note 23]{Weatherall}. 
 
 The investigations in the present paper are also relevant to the so-called syntax-semantics debate in philosophy of science. The issue here is, roughly, whether it is better to consider theories occurring in science as collections of linguistical objects (e.g., sentences of a given language), or as collections of structural objects of some kind. For a summary of the debate see Lutz \cite{Lutz} and Hudetz \cite{Hudetz}.  In this context, the need for a semantic characterization of definitional equivalence was raised in Halvorson~\cite{Halwhat}. Glymour~\cite{Gly} pointed out that de Bouv\`ere~\cite{Bou} contains such a characterization. Theorem~\ref{main-thm} in the present paper is another such semantic characterization. An advantage of Theorem~\ref{main-thm} is that it gives intuition about what properties of theories are preserved by definitional equivalence. Namely, by Theorem~\ref{main-thm}, a property of a theory is preserved when it can be expressed in terms of universes, isomorphisms and ultraproducts of models. Glymour \cite[p.296]{Gly} conjectures that each of the following four properties is preserved by classical definitional equivalence: having a one-element model, the model class being closed under substructures, the model class being closed under unions of chains, and having an equational axiomatization. Of these, the first property is clearly preserved by definitional equivalence because it is expressed by using the universes of the models. We show, after Theorem~\ref{postest-cor}, that neither one of the remaining three properties is preserved by classical definitional equivalence. 

Halvorson \cite[section 7]{Halwhat} proposes the programme to investigate what structure a model class naturally has and Glymour \cite[p.297]{Gly} appreciates this programme. This programme involves to endow the model class of a theory in such a way that from this structure on the model class, the theory can be recovered up to definitional equivalence. For propositional logic, Stone-duality provides such a structure in form of the Stone-topology on the model class. Stone duality has been generalized to first-order logic by several authors, e.g., Makkai~\cite{Makkai87} and Awodey and Forssell~\cite{Awo}. Halvorson points out the relevance of Stone duality for his programme and he mentions \cite{Makkai87} and \cite{Awo}.
Now, from the model-structures proposed in these two papers, the first-order theory can be recovered only up to the weaker many-sorted definitional equivalence. Theorem~\ref{main-thm} in the present paper suggests a structure on the model classes, we call this concrete ultracategory, from which a theory can be recovered up to classical definitional equivalence (and not only up to many-sorted definitional equivalence). See Remark 7.  We do not know of any other structure proposed in the literature on the model classes from which a theory can be recovered up to classical definitional equivalence. 

Example 2  points to an interesting difference between structural and language-based equivalences of theories.
Namely, Example 2 contains pairs of theories which are not equivalent with respect to any finitely-linguistic-based equivalence (see the proof of Lemma \ref{def-cor}), yet there is a bijection between their model classes that preserves isomorphisms and ultraproducts up to isomorphism. If such a bijection preserves ultraproducts not only up to isomorphism, then it establishes definitional equivalence according to Theorem \ref{main-thm}. This shows that preserving ultraproducts only up to isomorphism, which structural properties usually do, is not enough for establishing classical definitional equivalence.\bigskip

{\bf On the approach taken in the present paper.}
It is known that definability and automorphisms are intimately connected. Though it is not true that a relation is definable in a model if and only if all automorphisms of the model preserve the relation, something close is true: a relation is definable if and only if all automorphisms of all ultrapowers preserve (the corresponding ultrapower of) the relation (see \cite[Lemma 6.7.5]{pezs}). This theorem has proved to be quite useful so far for establishing definability and non-definability of relations. 

This paper can be viewed as a search for a similar complete method for establishing definitional equivalence and inequivalence of theories.  Section 2 contains two examples. The warm-up Example 1 shows that having same classes of automorphism groups does not entail definitional equivalence. It also motivates the notion of spectrum of concrete automorphism groups. Example 2 shows that having same spectrum of concrete automorphism groups still does not entail definitional equivalence. It also shows the importance of preserving ultraproducts. Section 3 contains a purely semantic characterization of definitional equivalence (Theorem \ref{main-thm}), which is also a complete method for establishing definitional equivalence by using concrete automorphism groups and ultraproducts. We then show how to use this method for establishing definitional inequivalence of two theories from Example 2 (Theorem \ref{postest-cor}). Finally, we make connections with related recent philosophy of science papers.

If not stated otherwise, we use the notation of Chang and Keisler \cite{CK}.

\section{Testing with automorphism groups}

We are in first-order logic. Two theories $\T_1$ and $\T_2$ are
said to be \emph{definitionally equivalent} when there are copies of
these theories with disjoint languages which have a joint
definitional extension. A copy of a theory $\T$ is a theory $\T'$
which is obtained from $\T$ by renaming some elements of the
vocabulary. A definitional extension of a theory is the theory where
some defined relations are added to the language. For discussion of this definition of 
definitional equivalence of theories see
the introduction and \cite[Definitions 10, 19, Theorem 4]{LSz}.
Two theories are said to be \emph{definitionally inequivalent} when they are not definitionally equivalent.
When $\T$ is a theory, $\Mod(\T)$ denotes the class of its models,
and when $\K$ is a class of similar models, $\Th(\K)$ denotes its
theory, i.e., the set of formulas valid in it. When $\Mm$ is a
model, $\Aut(\Mm)$ denotes its concrete automorphism group, i.e.,
the universe of $\Aut(\Mm)$ is the set of all \emph{automorphisms}
of $\Mm$ (i.e., permutations of the universe of $\Mm$ which leave
all relations of $\Mm$ unchanged as sets) and the sole operation of $\Aut(\Mm)$ is
the operation of composition.
\[ \Aut(\T) = \{\Aut(\Mm) : \Mm\in\Mod(\T)\}.\]

We begin with two examples. The first example serves to show that
searching for automorphism groups occurring in one but not the other
of the theories is not a complete method for showing failure of
definitional equivalence. 
\bigskip

\noindent {\bf Example 1.}(definitionally inequivalent theories
with same automorphism groups) We present theories $\T_1$ and
$\T_2$ such that $\Aut(\T_1)=\Aut(\T_2)$ and $\T_1$ is not
definitionally equivalent to $\T_2$. The two theories have the same
language, this language contains two binary relation symbols $\Es,
\R$.  The first theory, $\T_1$, states that at most one of $\Es$ and
$\R$ can be non-empty. The second theory, $\T_2$, states in addition
that when $\R$ is non-empty it is asymmetric:
\[ \T_1=\{ \forall xy\lnot\Es(xy)\lor\forall xy\lnot\R(xy)\}, \quad
\T_2=\T_1\cup\{\forall xy(\R(xy)\to\lnot\R(yx))\} .\]

The two theories have same automorphism groups because of the following. 
Let $\G$ denote
the class of automorphism groups of all models with one binary
relation, i.e., $\G=\{\Aut(\langle M,S\rangle) : S\subseteq M\times
M\}$. Clearly, $\Aut(\T_1)=\Aut(\T_2)=\G$ because in any model of
$\T_1$ or $\T_2$ there is at most one nonempty relation and the
empty relation does not affect the automorphism group, so
$\Aut(\T_1)\cup\Aut(\T_2)\subseteq\G$. The other containment follows
from the fact that neither of the theories make any restriction on
$\Es$.

To show that $\T_1$ and $\T_2$ are not definitionally equivalent, we
will exhibit a  concrete group  $\Gg$ that occurs as the
automorphism group for finitely many models altogether, but more
models of $\T_1$ than of $\T_2$ have $\Gg$ as their automorphism
group. 
Let the universe of $\Gg$ consist of one member, the identity
map on $H=\{ 0,1\}$. There are 12 binary relations on $H$ altogether
whose automorphism group consists only of the identity on $H$, 2 of
these are asymmetric. Thus there are 24 models in $\Mod(\T_1)$ with
automorphism group $\Gg$, because in each such model of $\T_1$ either $S$ is 
empty and $R$ is one of the 12 binary relations or the other way round. 
However, only 14 models in $\Mod(\T_2)$ has
$\Gg$ as automorphism group because either $R$ is empty and $S$ is one of 
the 12 above, or $S$ is empty and $R$ is one of the 2 antisymmetric relations.
This shows that there is no bijection between the models  of $\T_1$ and $\T_2$ 
which 
is such that corresponding models have the same 
automorphism group. Therefore, they 
are not model meargeable and so not definitionally equivalent. 

It may be interesting to have only infinite models for our theories.
An easy modification of $\T_1$ and $\T_2$ will do. Namely, we add
both to $\T_1$ and to $\T_2$ the infinitely many sentences that
together state that their models are infinite. We then have to
modify $\Gg$. The universe of the new $\Gg$ consists of
all permutations on $H=\{ 0,1,2,\dots\}$, the set of non-negative
integers, that leave 0 fixed. $\Box$
\bigskip

The previous example suggests that multiplicity of concrete automorphism
groups has to be taken into account when testing definitional
equivalence. We define the {\it spectrum of concrete automorphism groups} of a theory $\T$ as
a function that to each permutation group associates the number of non-isomorphic models of $\T$ that have this group as concrete automorphism group, i.e.,

\[ \AutSpec(\T)\de \{\langle \Gg, \nu(\Gg,\T)\rangle : \Gg \mbox{ is a permutation group }\}\]
where
\[\nu(\Gg,\T) \de \vert\{\Mm\in\Mod(\T) : \Aut(\Mm)=\Gg\}\slash\!\!\cong\vert.\]
Note that if two models have the same concrete automorphism group then they must have the same universe.

Definitionally equivalent theories have same spectrum of concrete automorphism groups. 
Therefore, for two theories to be definitional equivalent, it is necessary that they have same spectrum of concrete automorphism groups.
The most natural way of ensuring this is to require a bijection between their classes of models which preserves concrete automorphism groups as well as isomorphisms. This leads to the notion of a category of models formed from the models of a theory. 

The most common way of forming a \emph{category from the models} of
a first-order logic theory is to take the models of the theory as the
objects of the category and take the elementary embeddings%
\footnote{For the definition of  elementary embedding see \cite[p.84]{CK}.} between
these models as morphisms of the category. Let $\Cat(\T)$ denote
this category of models of $\T$. 
%
Often, it is useful to investigate a category of models with fewer morphisms taken into account. The \emph{model-iso-category} $\Cat^{iso}(\T)$ of a theory is defined by having $\Mod(\T)$ as its class of objects and having as morphisms only the isomorphisms between models.
The arguments in James Owen Weatherall \cite{Weatherall} point in
the direction to deal with the category of models when only
isomorphisms are taken as arrows, and not all elementary embeddings.
The idea is that in many realistic cases, just as ones dealt with in
\cite{Weatherall}, the scientific theory is not defined by a first-order logic theory, yet one has a clear sense of what models and
isomorphisms between these models can be.

Model categories come with a natural \emph{forgetful functor} to the category $\Set$ of all sets. These functors assign the universe $M$ to a model $\Mm$ and they assign the ``function content" to a morphism between two models. These are so natural in model theory that they are called \emph{the} forgetful functor. For definitions see \cite[Definition  5.1 (1)]{AHS}. A functor $F$ between model categories is called a \emph{concrete functor} iff it commutes with these natural forgetful functors. Thus a functor $F$ between model categories is a concrete one iff the universes of connected models are the same and if connected morphisms are the same as functions between the universes of models. Two model categories are called \emph{concretely isomorphic} iff there is a concrete isomophism between them.

Existence of concrete isomorphism between model-iso-categories is a natural generalization of having the same spectrum of concrete automorphism groups. The next theorem says that, in fact, it is not a generalization.

\begin{thm}
	Two theories have same spectrum of concrete automorphism groups if and only if their model-iso-categories are concretely isomorphic.
\end{thm}

\noindent{\bf Proof.} 
Let $\T_1$ and $\T_2$ be first-order theories and assume that $\AutSpec(\T_1)=\AutSpec(\T_2)$.  We are going to define a concrete isomorphism $b$ between their model-iso-categories.

The identity element of a permutation group is always of the form $\{(a,a) : a\in A\}$ for some $A$, let us call this $A$ the base of the permutation group. Let $\Gg,\Hh$ be permutation groups, let $h:A\to B$ be a bijection between the bases of $\Gg$ and $\Hh$, and define $\overline{h}(g) = h\circ g\circ h^{-1}$  for all $g\in G$. Then it is easy to see that $\overline{h}$ is an isomorphism between $\Gg$ and $\Hh$, we say that it is the base-isomorphism induced  by $h$. A \emph{base-isomorphism} between two permutation groups $\Gg,\Hh$ is an isomorphism betwen them that is induced by some $h$. We will also use the fact that if $h:\Mm\to\Nn$ is an isomorphism between the structures $\Mm,\Nn$, then $\overline{h}$ is a base-isomorphism between their automorphism groups.

Let $\cG$ be a class of representatives for the  base-isomorphism classes of permutation groups. That is, each permutation group has a base-isomorphic copy in $\cG$ and the elements of $\cG$ are pairwise non-base-isomorphic.
For any permutation group $\Gg\in\cG$ choose $\nu(\Gg,\T_1)$-many non-isomorphic models $\Mm(\Gg,i)$ of $\T_1$, for $i<\nu(\Gg,\T_1)$, and similarly choose $\nu(\Gg,\T_2)=\nu(\Gg,\T_1)$ non-isomorphic models $\Mm'(\Gg,i)$ of $\T_2$, with concrete automorphism group $\Gg$.  Then the models $\Mm(\Gg,i)$ for $\Gg\in\cG$ are pairwise non-isomorphic, i.e., $\Mm(\Gg,i)\cong\Mm(\Hh,j)$ for some $\Gg,\Hh,i,j$ implies $\Gg=\Hh$ and $i=j$. Similarly, the models $\Mm'(\Gg,i)$ are pairwise non-isomorphic.

Let $\Mm\in\Mod(\T_1)$. There is a unique $\Mm(\Gg,i)$ isomorphic to $\Mm$, as follows. Let $\Hh$ be the concrete automorphism group of $\Mm$ and let $\Gg\in\cG$ be base-isomorphic to $\Hh$ via the base-isomorphism $\overline{h}:\Hh\to\Gg$. Then the automorphism group of $h(\Mm)$ is $\Gg\in\cG$, thus $h(\Mm)$ is isomorphic to $\Mm(\Gg,i)$ for some $i$, by our construction.
Choose any isomorphism $f$ mapping $\Mm(\Gg,i)$ to $\Mm$ and let us define
\[ b(\Mm) \de f(\Mm'(\Gg,i)) .\]
We show that $b(\Mm)$ is well-defined, i.e., it does not depend on which isomorphism $f$ we choose.  Let $g$ be any other isomorphism between $\Mm(\Gg,i)$ and $\Mm$, we show that $f(\Mm'(\Gg,i))=g(\Mm'(\Gg,i))$. Indeed, $g=f\circ\alpha$ for $\alpha=f^{-1}\circ g\in\Aut(\Mm(\Gg,i))=\Gg$. But $\alpha\in\Gg=\Aut(\Mm'(\Gg,i))$, so $g(\Mm'(\Gg,i))=f(\alpha(\Mm'(\Gg,i)))=f(\Mm'(\Gg,i))$. 

We define $b$ on the morphisms. Let $h:\Mm\to\Nn$ be an isomorphism between $\Mm,\Nn\in\Mod(\T_1)$. We have seen that $f:\Mm(\Gg,i)\to\Mm$ for some $f,\Gg,i$ and so $g:\Mm(\Gg,i)\to\Nn$ for $g=h\circ f$. Thus, by definition, $b(\Mm)=f(\Mm'(\Gg,i))$ and $b(\Nn)=g(\Mm'(\Gg,i))$. Hence, $h:b(\Mm)\to b(\Nn)$ is an isomorphism by $g\circ f^{-1}=h\circ f\circ f^{-1}$. We define
\[ b(h) \de h .\] 

We now show that $b$ is an isomorphism between the model-iso-categories of $\T_1$ and $\T_2$.
First we show that the function $b:\Mod(\T_1)\to\Mod(\T_2)$ defined this way is a bijection between $\Mod(\T_1)$ and $\Mod(\T_2)$. Indeed, let $\Mm'\in\Mod(\T_2)$ be any model. There is a unique $\Mm'(\Gg,i)$ isomorphic to it, say via $f:\Mm'(\Gg,i)\to\Mm'$. Let $\Mm=f(\Mm(\Gg,i))$, then $\Mm'=b(\Mm)$, by the definition of $b$. Thus, the range of $b$ is $\Mod(\T_2)$. To see that $b$ is one-to-one, let $\Mm,\Nn\in\Mod(\T_1)$. Assume that $b(\Mm)=b(\Nn)$. By the definition of $b$, there are $\Mm(\Gg,i), \Mm(\Hh,j)$ and isomorphisms $f:\Mm(\Gg,i)\to \Mm$, $g:\Mm(\Hh,j)\to\Nn$ such that  $b(\Mm)=f(\Mm'(\Gg,i))$ and $b(\Nn)=g(\Mm'(\Hh,j))$. By $b(\Mm)=b(\Nn)$ then $\Mm'(\Gg,i)$ is isomorphic to $\Mm'(\Hh,j)$,  therefore, $(\Gg,i)=(\Hh,j)$ and  $f(\Mm'(\Gg,i))=g(\Mm'(\Gg,i))$. Thus $f^{-1}\circ g\in\Aut(\Mm'(\Gg,i))=\Gg$. So,  $\Mm=f(\Mm(\Gg,i))=f((f^{-1}\circ g)(\Mm(\Gg,i)))=g(\Mm(\Gg,i))=g(\Mm(\Hh,j))=\Nn$. 

We turn to the proof for $b$ being a bijection between the set of isomorphisms from $\Mm$ to $\Nn$ and the set of isomorphisms from $b(\Mm)$ to $b(\Nn)$, for any $\Mm,\Nn\in\Mod(\T_1)$. To show surjectivity, let $h:b(\Mm)\to b(\Nn)$. By the definition of $b(\Mm)$, we have that $\Mm=f(\Mm(\Gg,i))$ and $b(\Mm)=f(\Mm'(\Gg,i))$, for some $f,\Gg,i$. Thus, $f:\Mm'(\Gg,i)\to b(\Mm)$, and so $h\circ f:\Mm'(\Gg,i)\to b(\Nn)$, by $h:b(\Mm)\to b(\Nn)$. Let $\Nn'=f(h(\Mm(\Gg,i))$. By the definition of $b(\Nn')$ then  $b(\Nn')=(h\circ f)(\Mm'(\Gg,i))=b(\Nn)$. Thus $\Nn'=\Nn$ because $b$ is one-to-one on $\Mod(\T_1)$, i.e., $\Nn=(h\circ f)(\Mm(\Gg,i))=h(f(\Mm(\Gg,i)))=h(\Mm)$. Thus, $h:\Mm\to\Nn$ is an isomorphism and $b(h)=h$. By definition, it is  clear that $b$ is one-to-one on the morphisms, and also that it preserves composition of morphisms both directions. This finishes the proof for $b$ being a category theoretical isomorphism between the model-iso-categories of $\T_1$ and $\T_2$. It is concrete, by its definition.

In the other direction, assume that $b$ is a concrete isomorphism between $\Cat^{iso}(\T_1)$ and $\Cat^{iso}(\T_2)$. Then $\Aut(\Mm)=\Aut(b(\Mm))$, and $\Mm\cong\Nn$ iff $b(\Mm)\cong b(\Nn)$, for all $\Mm,\Nn\in\Mod(\T_1)$. Therefore, $\AutSpec(\T_1)=\AutSpec(\T_2)$. 
\qed\bigskip

The next example shows that having same spectrum of automorphism groups still does not entail definitional
equivalence.
It is more refined than the previous one. We will see
that it shows, in a sense, a limit till we still can get failure of
definitional equivalence (compare Lemma \ref{functor-lem} with Theorem~\ref{main-thm}).
It also serves as a counterexample to Barrett and Halvorson's
conjecture that, among first-order logic theories with finite
signatures, categorical equivalence implies many-sorted (Morita)
definitional equivalence. With this, we answer in the negative
\cite[Question 6.1]{BH16} as well as \cite[Question 1,
p.77]{BDis17}.
\bigskip

\noindent {\bf Example 2.}(uncountably many theories with same model category) We present continuum many complete theories 
on a finite similarity type with same automorphism spectrum such that no 
two of them are definitionally equivalent. 
Moreover, their model categories are isomorphic via concrete functors
which preserve ultraproducts up to
isomorphism, and further, no two of the theories are even many-sorted definitionally (Morita) equivalent. 
(The latter notion will be introduced later, below Lemma \ref{functor-lem}.)

We work in the similarity type which contains one constant symbol
$\zero$, one unary function symbol $\suc$, and one unary relation
symbol $\R$.  
Let $n$ be a natural number, then $\suc^n(x)$
denotes the term where $\suc$ is $n$-times applied to $x$, i.e.,
$\suc^0(x)=x$ and $\suc^{(n+1)}(x)=\suc(\suc^n(x))$. 
For each subset  $S$ of the natural numbers $\omega$ let
\[ \T(S)\de\{ \R(\suc^n(\zero)) : n\in S\}\cup \{\lnot\R(\suc^n(\zero)) : n\notin S\}\cup\Th(\langle\omega,\zero,\suc\rangle)\]
where $\langle\omega,\zero,\suc\rangle$ denotes natural numbers $\omega$ with zero as $\zero$ and the successor function as $\suc$. 

A set $S$ of natural numbers is called \emph{irregular} if all finite patterns occur in it. In more detail, let $n>0$ be a positive number and let $P\subseteq\{0,1,...,n-1\}$. We say that the $P,n$-pattern occurs at $x$ in $S$ if $\{ m<n: suc^m(x)\in S\}=P$.
For example, $S=\{ 0,2,4,6,...\}$ is not irregular, because the pattern $\{0,1\},2$ does not occur in it (i.e., $x,\suc(x)\in S$ does not hold for any $x\in\omega$). 

There are continuum many irregular subsets of $\omega$. This can be seen as follows.
Construct an infinite sequence of $0,1,x$ by first laying the two $0,1$-sequences of length 1 after each other in alphabetical order, then mark the next number by an $x$, then lay the four $0,1$-sequences of length 2 after each other in alphabetical order and mark the next number by an $x$, etc. This sequence will begin like $\langle 0,1,x,0,0,0,1,1,0,1,1,x, 0,0,0,...\rangle$. There are infinitely many $x$s in this sequence and so there are continuum many ways of replacing the $x$s with $0$ or $1$.  Each of the continuum many $0,1$-sequences that are obtained this way is a characteristic function of an irregular set. This proves that there are at least continuum many irregular sets. There can be at most continuum many irregular subsets of $\omega$ since there are continuum many subsets of $\omega$.

We are going to show that the model categories $\Cat(\T(S))$
for irregular sets $S$ are isomorphic to each other in a strong constructive way, see Lemma \ref{functor-lem}.
\bigskip

We say that $\Nn$ is an \emph{induced subalgebra} of $\Mm$ when the
$\R$-free part of $\Nn$ is a subalgebra of the $\R$-free part of
$\Mm$ and the $\R$-relation of $\Nn$ is that of $\Mm$ restricted to
the universe of $\Nn$. For the definition of elementary submodel see \cite[p.84]{CK}.

\begin{lem}\label{irreg-lem} Let $S\subseteq\omega$ be irregular. Then (i)-(ii) below hold.
\begin{description}
\item{(i)}
The elementary submodels of a model of $\T(S)$ are exactly its induced
subalgebras.
\item{(ii)}
$\T(S)$ is a complete theory.
\end{description}
\end{lem}

\noindent {\bf Proof.} 
	Let $\N$ denote the set of natural numbers
with $0$ as constant $\zero$ and the successor function as unary
distinguished function $\suc$, and let $\Z$ denote the set of integers with
the successor function as unary distinguished function  $\suc$.  Note that  $\Z$ does not have $\zero$ in its language. Any model of $\Th(\langle\omega,\zero,\suc\rangle)$  is a disjoint union
of one copy of $\N$ together with some copies of $\Z$. 
 When $k$ is negative, $\suc^{k}(x)=y$ means $\suc^{-k}(y)=x$, we say that $\suc^{k}(x)$ exists when such a $y$ exists. In models of $\Th(\langle\omega,\zero,\suc\rangle)$ such a $y$ is unique when it exists.  When $\Nn$ is a model of $\Th(\langle\omega,\zero,\suc\rangle)$, by a \emph{$\Z$-part} of $\Nn$ we mean a subset of $N$ of the form $\{\suc^n(a) : n\in\omega\}\cup\{\suc^{-n}(a) : n\in\omega\}$ for some $a\in N$.
 By a
$\Z$-model we mean $\Z$ together with a unary relation $R$ and by $\langle\N,S\rangle$ 
we mean $\N$ expanded with $S$ as the unary relation $\R$. We are going to
prove the following statement (*). 

In (*) as well as later on, we will use ultraproducts (\cite[ch.4]{CK}).
As in \cite{CK}, when $U$ is an ultrafilter on the set $I$ and $\langle\Mm_i : i\in
I\rangle$ is an $I$-sequence of similar models,
$\prod_{U}\langle\Mm_i : i\in I\rangle$, or sloppily just
$\prod_{U}\Mm_i$, denotes the
$U$-ultraproduct of the models $\Mm_i$ and $y_U$ denotes the equivalence-class of $y$ in $\Pi_U\Mm_i$, for $y\in\Pi_{i\in I}M_i$.
When each $\Mm_i=\Aa$ for some $\Aa$, we call $\prod_{U}\Mm_i$ an
ultrapower of $\Aa$ and we denote it by $\Pi_U\Aa$.

\begin{description}
\item{(*)}
Assume that $\Mm$ is a countable model of $\T(S)$ and $U$
is a nonprincipal ultrafilter on a countable set $I$. Then
$\Pi_U\Mm$ is isomorphic to a disjoint union of a copy of $\langle\N,S\rangle$ with continuum
many copies of each possible $\Z$-model.
\end{description}
Indeed, (*) is true because each $\Z$-model can be put together in
the ultrapower from its finite parts which are patterns occurring in
$\Mm$, and in fact, each such pattern occurs infinitely many
times in $\Mm$. In more detail:
Let $\langle\Z,R\rangle$ be any $\Z$-model, we show that continuum 
many disjoint copies of it occurs in the ultrapower of $\Mm$. We may assume that $I=\omega$ because $I$ is countable. For each $n>0$ let $R_n\de\{m\le2n : m-n\in R\}$. The pattern $R_n,2n+1$ occurs in $S$ because $S$ is irregular. In fact, each pattern occurs in an irregular set infinitely many times because each finite pattern has infinitely many different extensions to other finite patterns and each of these patterns occur in the irregular set. Let $X_n$ be the set of $x$s where $R_n,2n+1$ occurs in $S$ and let $Y_n\de\{x+n : x\in X_n\}$. First we show that in $\Pi_U\Mm$ each element of $\Pi_U Y_n$ lies on a copy of $\langle\Z,R\rangle$.  Indeed, let  $x_n\in X_n$ and $y_n\de x_n+n$ for all $n\in\omega$. Let $y\de\langle y_n : n\in\omega\rangle$, and let $k\in\Z$ be arbitrary. We will show that $\suc^k(y_U)$ exists and $k\in R$ iff $R(\suc^k(y_U))$ in $\Pi_U\Mm$. By our definitions, for all $n$ such that $2n\ge k$ we have that $k\in R$ iff $k+n\in R_n$ iff $x_n+k+n\in S$ iff $y_n+k\in S$ iff $R(\suc^k(y_n))$ in $\Mm$. Since $U$ is nonprincipal on $I=\omega$, this means that $R(\suc^k(y_U))$ in $\Pi_U\Mm$. We have seen that $y_U$ is in a copy of $\langle\Z,R\rangle$ for all $y\in\Pi_{n\in\omega} Y_n$. Since each $Y_n$ is countably infinite, the cardinality of $\Pi_U Y_n$ is continuum (see, \cite[Prop.4.3.9]{CK}). Since each copy of $\langle\Z,R\rangle$ is countable, this means that $\Pi_U\Mm$ contains continuum many disjoint copies of $\langle\Z,R\rangle$, and we are done with proving (*).

Proof of (i): An  elementary submodel of $\Mm$ has to be an
induced subalgebra. Conversely, assume that $\Nn$ is an induced
subalgebra of $\Mm$, we show that it is an elementary submodel.
We will use the testing method in \cite[Prop.3.1.2]{CK}. Thus,
assume that $\varphi(\bar{x},y)$ is a first-order logic formula in
the language of $\Mm$, assume that $\bar{a}$ is an appropriate
sequence of elements of $\Nn$, and $\Mm\models\exists
y\varphi(\bar{a},y)$. We have to show the existence of $a'\in\Nn$
such that $\Mm\models\varphi(\bar{a},a')$. We have $\Pi_U\Mm\models\exists y\varphi(d(\bar{a}),y)$ since the diagonal (or
natural) embedding $d$ of a model into its ultrapower is an
elementary one \cite[Cor.4.1.13]{CK}. Let $b\in\Pi_U\Mm$ be
such that $\Pi_U\Mm\models\varphi(d(\bar{a}),b)$.
Now, $\Pi_U\Nn$ is an induced subalgebra of $\Pi_U\Mm$, by
$\Nn$ being an induced subalgebra of $\Mm$. There are infinitely
many $\Z$-parts in $\Pi_U\Nn$ that do not contain any element
of $d(\bar{a})$ and that are isomorphic to the $\Z$-part of
$\Pi_U\Mm$ containing $b$, by (*). Take an automorphism of
$\Pi_U\Mm$ that interchanges the $\Z$-part of $b$ with any of
such a $\Z$-part of $\Pi_U\Nn$ and leaves anything else fixed.
There is such an automorphism by the choice of the $\Z$-part of
$\Pi_U\Nn$ and since $\Mm\in\Mod\T(S)$. Let $c$ be the
image of $b$ under such an automorphism, then $\Pi_U\Mm\models\varphi(d(\bar{a}),c)$, since the automorphism leaves the
elements of $d(\bar{a})$ fixed. Then $\Mm\models\varphi(\bar{a},a')$
for some $a'\in\Nn$ by the fundamental theorem of ultraproducts
\cite[Thm.4.1.9(ii)]{CK} since $c\in\Pi_U\Nn$. We have shown
that $\Nn$ is an elementary submodel of $\Mm$. 

Proof of (ii): Assume that $\Mm,\Nn\in\Mod\T(S)$, we have to show that $\Nn$ is elementarily
equivalent to $\Mm$. We may assume that $\Mm$ and $\Nn$ are
countable, by the downward L\"owenheim--Skolem--Tarski theorem
(\cite[Cor.2.1.4]{CK}). Now, $\Mm$ and $\Nn$ are elementarily
equivalent by (*), since they have isomorphic ultrapowers.
 The proof of
Lemma~\ref{irreg-lem} is complete. 
\bigskip

By using Lemma~\ref{irreg-lem}, we now specify a functor $F$ between
the model categories of $\T(S)$ and $\T(Z)$, for any irregular sets $S$ and $Z$.
Let $\Mm=\langle M, \zero, \suc, R\rangle\in\Mod(\T(S))$.
We define
\[ F(\Mm)\de\langle M, \zero, \suc, (R\setminus\{\suc^{k}(\zero) : k\in S\})\cup\{\suc^{k}(\zero) : k\in Z\}\rangle .\]
That is, $F(\Mm)$ is defined to be $\Mm$ except that $R$ on the $\N$-part of $\Mm$ is
changed to be the $R$ of the $\N$-part of a $\T(Z)$ model. For an elementary embedding $f:\Mm\to\Nn$ between $\Mm,\Nn\in\Mod(\T(S))$ let us define
\[ F(f)\de f.\]

\begin{lem}\label{functor-lem} Let $S$ and $Z$ be irregular sets and let $F$ be the function defined above.
	\begin{description}
	\item[(i)] $F$ is a concrete isomorphism between $\Cat(\T(S))$ and $\Cat(\T(Z))$.
	\item[(ii)] $F$ preserves ultraproducts of models up to isomorphism, i.e., $F$ takes an ultraproduct of models of $\T(S)$ to a model isomorphic to the corresponding ultraproduct of the $F$-images of the models.	
	\end{description}
\end{lem}

\noindent{\bf Proof.} 
 $F$ is a functor, since ($f$ is an elementary embedding of $\Mm$ into $\Nn$ if and only if it is an elementary embedding of $F(\Mm)$ into $F(\Nn)$),  by Lemma~\ref{irreg-lem} and the construction of $F$. Thus $F$ is a concrete isomorphism by its construction.

To show that $F$ preserves ultraproducts up to isomorphism, let $U$ be an ultrafilter on a set $I$ and let $\Mm_i\in\Mod(\T(S))$ for all $i\in I$. We will define an isomorphism  $j$ between $F(\Pi_U\Mm_i)$ and $\Pi_U F(\Mm_i)$. Let  $\Nn_i$ denote the $\N$-part of $\Mm_i$, for each $i\in I$. Then each $\Nn_i$ is isomorphic to $\langle \N,\{\suc^k(\zero) : k\in S\}\rangle$ by $\Mm_i\in\Mod\T(S)$.  Let $y\de\langle y_i : i\in I\rangle\in\Pi_{i\in I}M_i$. We define
\[ j(y_U) \de y_U\quad\mbox{if}\quad \{ i\in I : y_i\notin N_i\}\in U .\]
To define $j$ on the rest, assume first that $U$ is not $\omega^+$-complete. Then by a straightforward modification of the proof of (*) we get that both $\Pi_U\Nn_i$ and $\Pi_U F(\Nn_i)$ consist of one $\N$-model together with continuum many copies of all possible $\Z$-models. If $U$ is $\omega^+$-complete, then both $\Pi_U\Nn_i$ and $\Pi_U F(\Nn_i)$ consist of one $\N$-model only by \cite[Prop.4.2.4]{CK}. 
In both cases there is an isomorphism between $F(\Pi_U\Nn_i)$ and $\Pi_U F(\Nn_i)$. We define
\[ j\mbox{ be any isomorphism between }F(\Pi_U\Nn_i)\mbox{ and }\Pi_U F(\Nn_i) \]
and be identity on the rest. It is not difficult to check that $j:F(\Pi_U\Mm_i)\to\Pi_U F(\Mm_i)$ is an isomorphism. 
This finishes the proof of Lemma \ref{functor-lem}. 
\bigskip

We have seen that, for any two irregular sets $S$ and $Z$, the model categories of $\T(S)$ and $\T(Z)$ are rather close to each other in a constructive way. We now turn to definability issues between $\T(S)$ and $\T(Z)$. In logic, there are two weaker versions of definitional equivalence between theories in use. One is called many-dimensional (\cite{Visser, harnik}) or many-sorted (\cite{Mad02, AN}) definitional equivalence, and it is also called Morita equivalence of theories (\cite{BH16, Halvorson}). The other is called bi-interpretability between theories 
(\cite{Hodges, Visser}). Both notions are weaker than definitional equivalence between first-order logic theories in the sense that when $\T_1$ and $\T_2$ are definitionally equivalent then they are also many-dimensionally equivalent and bi-interpretable. For a comparison of these notions see Barrett and Halvorson \cite{BH16}. We will rely on the deifinitions in the mentioned references, we do not recall them.

\begin{cor}\label{def-cor}\mbox{}
	\begin{description}
		\item[(i)] All the theories $\T(S)$ with $S$ irregular have same spectrum of automorphism groups.
		\item[(ii)] There is an uncountable set $\mathcal S$ of irrregular sets such that no $\T(S)$ and $\T(Z)$ for distinct $S,Z\in\mathcal S$ are definitionally equivalent, many-sorted definitionally equivalent or bi-interpretable. 
	\end{description}	
\end{cor}

\noindent{\bf Proof.}
	 (i) follows from Lemmas \ref{irreg-lem} and \ref{functor-lem}. Each of definitional equivalence, many-sorted definitional equivalence and bi-interpretability of two theories can be specified by the use of finitely many formulas on the language of the theories, see the references given for their definitions. Therefore, a concrete theory can be definitionally equivalent to at most countably many theories on a given other similarity type. This implies that of the continuum many theories $\T(S)$ on the same language, there are continuum many pairwise non-equivalent theories (neither many-sorted equivalent nor bi-interpretable). This finishes the proof of Corollary~\ref{def-cor}. 
With this, presentation of Example 2 is finished. \quad $\Box$

\bigskip
The essence of Example 2 above is that the model categories of
$\T(S)$ for irregular sets $S$ are almost
the same because the $R$ on the $\N$-parts do not play a role in
this category. However, the $R$ on the $\N$-part can code more
``information" than available (syntactical) translations between theories and therefore many such theories have to be definitionally
inequivalent.
\bigskip

\noindent {\bf Remark 3.}($F$ does not preserve ultraproducts)  The functor $F$ constructed above Lemma \ref{functor-lem} does not preserve ultraproducts, it preserves ultraproducts only up to isomorphism. This follows from Theorem~\ref{main-thm} in the next section and Corollary~\ref{def-cor}(ii). We now want to provide a concrete example that shows that $F$ does not preserve ultraproducts. Recall the continuum many irregular sets constructed above Lemma~\ref{irreg-lem}. Let $S_0$ and $S_1$ be the irregular sets we obtain by filling all the $x$s with $0$ and and by filling all the $x$s with $1$, respectively. Then $S_0\subseteq S_1$ and $S_1\setminus S_0$ is infinite. Let $\Nn_i\de\langle\omega,\zero,\suc,R_i\rangle$ where $R_i=\{\suc^k(\zero) : k\in S_i\}$  for $i=0,1$.  Consider the functor $F$ between $\T(S_0)$ and $\T(S_1)$. Then $F(\Nn_0)=\Nn_1$ by definition of $F$. Let $X\subseteq\omega$ be an infinite set which is disjoint from $S_0$ but is contained in $S_1$, let $U$ be a nonprincipal ultrafilter on $I=\omega$ such that $X\in U$ and let $y=\langle\suc^k(\zero) : k\in\omega\rangle$. Then $R(y_U)$ does not hold in $F(\Pi_U\Nn_0)$ while $R(y_U)$ holds in $\Pi_U F(\Nn_0)$ showing that the two structures are not the same (though, isomorphic). We will see in the next section that in fact $\T(S_0)$ is not definitionally equivalent to $\T(S_1)$ because there is no concrete isomorphism between their model categories that would preserve ultraproducts, see Theorem \ref{postest-cor}.  $\Box$\bigskip

\noindent {\bf Remark 4.}(more striking example) We can modify the above example to give a
more striking counterexample to the conjecture in \cite{BH16} which
at the same time is analogous to the example in the proof of
\cite[Theorem 5.7]{BH16}. The similarity type of $\T_1$ and $\T_2$
will be as in Example 2. The first theory, $\T_1$ states only that
$\zero$ is not in relation $\R$:
 \[ \T_1 = \{\lnot\R(\zero)\} .\]
For defining $\T_2$, take any irregular set $S$ such that $0\in S$, and then $\T_2$ is
\[
\T_2 = \{ \R(\zero)\to\varphi\ \ :\ \ \varphi\in\T(S) \}.\]
That is, the models of $\T_2$ are those of $\T_1$ together with all
the models of $\T(S)$.
Now, $\T_1$ is finitely axiomatized while it is easy to see that
$\T_2$ cannot be axiomatized finitely (e.g., by showing that the
complement of $\Mod(\T_2)$ is not closed under ultraproducts). Since
intertranslatability is an essence of definitional equivalence both
for the classical and the many-sorted versions, as e.g., Halvorson
\cite{Halvorson} argues, being finitely axiomatized is preserved,
for theories of finite similarity types, by the weaker many-sorted
(Morita) definitional equivalence also. So, $\T_1$ and $\T_2$ are
not Morita definitionally equivalent. However, their model
categories are equivalent, in fact isomorphic, as in \cite[Theroem
5.7]{BH16}.: a model category consists of isolated islands of
$\Cat(\Th(\Mm))$ for the models $\Mm$ of the theory (because if
there is a morphism between $\Mm$ and $\Nn$ then $\Mm$ and $\Nn$ are
elementarily equivalent since this morphism is an elementary
embedding of $\Mm$ into $\Nn$). Now, by Lemma~\ref{irreg-lem}, the
extra island of $\Cat(\T_2)$ is isomorphic to any one of the
continuum many islands $\Cat(\T(Z))$ of $\Cat(\T_1)$ where $Z$
is an irregular set with $0\notin Z$.
$\Box$\bigskip

\section{Testing with automorphism groups and ultraproducts}

We are ready to turn to the positive results of this paper. Lemma~\ref{functor-lem} suggests that, besides automorphism groups, ultraproducts have to be taken into account in testing definitional equivalence. Indeed, Theorem~\ref{main-thm} below gives such a characterization making our search for a complete testing method successful.

The following theorem is a semantic characterization of definitional
equivalence. It is a slight modification of the Theorem in
\cite{vBP} which is a semantic characterization of restricted
interpretations between theories. For a closely related theorem see also Kochen \cite[Theorem 12.1]{Kochen}.

\begin{thm}\label{main-thm} Two theories $\T_1$ and $\T_2$ are definitionally
equivalent if and only if there is a bijection $b$ between their
model classes that satisfies the following two conditions.
\begin{description}
\item{(i)} An isomorphism between different models of $\T_1$ is
an isomorphism between their $b$-images and vice versa. In
particular, the universes of $\Mm$ and $b(\Mm)$ are the same.
\item{(ii)} Ultraproducts are preserved by $b$ in the sense that
$b(\prod_{U} \Mm_i) = \prod_{U} b(\Mm_i)$ for all ultrafilters $U$
and models $\Mm_i$ in $\Mod(\T_1)$.
\end{description}
\end{thm}

\noindent {\bf Proof.} 
	The proof follows that of
\cite[Theorem]{vBP}. Let assume first that the languages of $\T_1$
and $\T_2$ are disjoint. Assume that we have a bijection $b$
satisfying (i)-(ii). We define a class $\K$ of models in the
similarity type as the union of the similarity types of $\T_1$ and
$\T_2$ and we will show that the first-order logic theory of $\K$ is a
joint definitional extension for both $\T_1$ and $\T_2$. For a model
$\Mm\in\Mod(\T_1)$ let \[ \overline{\Mm}=\langle \Mm,
b(\Mm)\rangle\]  denote the model whose universe is the joint
universe of $\Mm$ and $b(\Mm)$, the relation and function symbols of
the language of $\T_1$ are interpreted as in $\Mm$, and the relation
and function symbols of the language of $\T_2$ are interpreted as in
$b(\Mm)$. Let
\[ \K = \{ \langle\Mm,b(\Mm)\rangle : \Mm\in\Mod(\T_1) \}  .\]
We will show that $\K$ is axiomatizable, i.e., $\K=\Mod\Th\K$. We
use \cite[Cor.6.1.16(i)]{CK} which states that a class is elementary
if and only if it is closed under taking ultraproducts and
isomorphic images, and the complement is closed under ultrapowers.
Now, $\K$ is closed under ultraproducts and isomorphisms by
conditions (ii) and (i), since $\Mod(\T_1)$ is elementary. Assume
that  $\Aa=\langle \Mm,\Nn\rangle$ is such that an ultrapower
$\Pi_U\Aa$ is in $\K$. We have to show that $\Aa\in\K$. Now,
$\Pi_U\Aa=\langle\Pi_U \Mm,\Pi_U\Nn\rangle$, and then
$\Pi_U\Aa\in\K$ means that $\Pi_U \Nn=b(\Pi_U \Mm)$. By
condition (ii) we have $b(\Pi_U\Mm)=\Pi_U b(\Mm)$. Thus
we have $\Pi_U\Nn=\Pi_U b(\Mm)$. This implies
$\Nn=b(\Mm)$ since any structure $\Bb$ can be recovered from
$\Pi_U\Bb$. We have seen that $\K$ is an elementary class, let
\[ \T = \Th(\K) .\]
Now, we show that $\T$ is a definitional extension of $\T_1$. When
the language of $\T_2$ has only one non-logical symbol, this follows
immediately from Beth's definability theorem (see
\cite[Thm.2.2.22]{CK}), since for each $\Mm\in\Mod(\T_1)$ there is
at most one relation satisfying $\T$, namely that of $b(\Mm)$.
However, a generalized version of Beth's theorem is well-known as
folklore: if the $\RR$-free reduct of each model of $\T$ can be 
extended to at most one model of $\T$, then $\T$ explicitly defines 
each member of $\RR$ by a formula on the language of the $\RR$-free reducts.%
\footnote{We give a short proof of this in the Appendix.}
The proof that $\T$ is a
definitional extension of $\T_2$ is completely analogous. Thus,
$\T_1$ and $\T_2$ are definitionally equivalent theories.

Assume now that the languages of $\T_1$ and $\T_2$ are not disjoint.
Rename the symbols in the language of $\T_2$ so that the new symbols
be distinct from any one used in $\T_1$ and $\T_2$, call the new
theory $\T_2'$. Now, there is a natural bijection
$b_1:\Mod(\T_2)\to\Mod(\T_2')$ satisfying conditions (i)-(ii), and
$b_1\circ b:\Mod(\T_1)\to\Mod(\T_2')$ also satisfies (i)-(ii). These
bijections are between models of theories of disjoint languages.
Apply the previous case to $b_1$ and $b_1\circ b$, and use that
definitional equivalence is a transitive relation by \cite{LSz}.
\qed\bigskip

\noindent{\bf Remark 5.}(relationship of Theorem~\ref{main-thm} with the van Benthem and Pearce result)
The theorem in \cite{vBP}, call it BP-theorem for van Benthem and Pearce theorem,  seems to be neither stronger nor weaker than Theorem~\ref{main-thm} above. It is not weaker because the kind of interpretation it deals with is restricted interpretation which is in between classical and Morita-interpretation. It is not stronger because it deals with interpretation and not with equivalence. In more detail, assume that there is a bijection between $\Mod(\T_1)$ and $\Mod(\T_2)$ satisfying (i),(ii) of Theorem~\ref{main-thm}.  By applying the BP-theorem, we get that there are two restricted interpretations, one from $\T_1$ to $\T_2$ and the other from $\T_2$ to $\T_1$.  However, we know that mutual interpretability even with strong properties does not imply definitional equivalence, see e.g., \cite{AMNMLQ}.  Although the BP-theorem does not seem to imply Theorem 2, the proof of Theorem 2 here is just a slight modification of the proof of the BP-theorem in \cite[p.296]{vBP}. \quad $\Box$
\bigskip

\noindent{\bf Remark 6.}(automorphism groups and elementary embeddings in Theorem \ref{main-thm})
 The word ``different" can be omitted from condition (i) of
Theorem~\ref{main-thm} and the theorem remains true. This is true
because condition (i) implies that the automorphism groups are
preserved by $b$ in the sense that $\Aut(\Mm)=\Aut(b(\Mm))$ for all
$\Mm\in\Mod(\T_1)$.
Indeed, if $\alpha\in\Aut(\Mm)$, then let $f:\Mm\to\Mm'$ be any
isomorphism where $\Mm'$ is different from $\Mm$, there is always
such an $f$. Then both $f$ and $\alpha\circ f$ are isomorphisms
between $b(\Mm)$ and $b(\Mm')$ by condition (ii), thus
$\alpha=\alpha\circ f\circ f^{-1}$ is an automorphism of $b(\Mm)$.
In a sense, this corollary about the automorphism groups is the
essential part of condition (i).

Also, Theorem \ref{main-thm} remains true if in (i) we require to preserve all 
elementary embeddings in place of all isomorphisms. The reason is that elementary 
embeddings are preserved by definitional equivalence.\quad $\Box$
\bigskip

The proof of the following theorem intends to illustrate the use of Theorem 
\ref{main-thm} for proving definitional inequivalence. Recall the definitions of $S_0$ and $S_1$ from Remark 3.

\begin{thm}\label{postest-cor}
	$\T(S_0)$ and $\T(S_1)$ are not definitionally equivalent.
\end{thm}

\noindent{\bf Proof.} 
	Let $b:\Mod(S_0)\to\Mod(S_1)$ be any bijection that preserves isomorphisms between distinct models. (We note that there is such a function $b$, see Lemma \ref{functor-lem}(i).) It preserves automorphism groups also, see Remark 3. We will show that $b$ cannot preserve all ultrapowers. By Theorem \ref{main-thm}, this will prove that $\T(S_0)$ and $\T(S_1)$ are not definitionally equivalent. 

Let $\Mm=\langle\omega,\zero,\suc,S_0\rangle\in\Mod(\T(S_0))$.  Let $b(\Mm)=\langle\omega,o,F,R\rangle$ and let $U$ be any nonprincipal ultrafilter on $\omega$. First we show that $b(\Pi_U\Mm)\ne\Pi_U b(\Mm)$ if $b(\Mm)$ contains any copy of a $\Z$-model. Indeed, assume that $F^k(n)$ exists for all $k\in\Z$ for some $n\in\omega$. Let $\langle\Z,P\rangle$ be the $\Z$-model that is isomorphic to the induced subalgebra of $b(\Mm)$ with universe $\{ F^k(n) : k\in\Z\}$.  By (*) in the proof of Lemma \ref{irreg-lem}, $\Pi_U b(\Mm)$ contains infinitely many copies of this $\Z$-model. Therefore, the image of the $\Z$-model in $b(\Mm)$ under the diagonal embedding can be interchanged with a distinct copy of this $\Z$-model in $\Pi_U b(\Mm)$. On the other hand, all automorphisms of $\Pi_U\Mm$ leave the diagonal embedding of $\Mm$ unchanged.  Thus $b(\Pi_U\Mm)$ cannot be $\Pi_Ub(\Mm)$ since the two have different automorphism groups. Therefore, we assume in the rest
\begin{description}\label{1}
\item{(1)} $\omega=\{ F^n(o) : n\in\omega\}$  and thus $R=\{ F^n(o) : n\in S_1\}$.	
\end{description}
\comment{Next we show that $b(\Pi_U\Mm)\ne\Pi_U b(\Mm)$ if there are infinitely many jumps in $F$, i.e., if 
\[J\de\{ k\in\Z : F(y)=\suc^k(y)\mbox{ for some }y\in\omega\}\]
 is infinite. Indeed, for all $k\in J$ let $y_k$ be such that $F(y_k)=\suc^k(y_k)$. Let $f:J\to\omega$ be any bijection (there is such because $J$ is countably infinite) and let $y\de\langle y_{k} : f(k)\in\omega\rangle$. Then $y$ is an $\omega$-sequence, namely $y:\omega\to\omega$, so $y_U\in\Pi_U\omega$. Further, $F(y_U)$ is not in the copy of the $\Z$-model containing $y_U$ in $\Pi_U\Mm$,  i.e., $F(y_U)\ne\suc^k(y_U)$ for all $k\in\Z$, by the definition of $J$. Choose an automorphism of $\Pi_U\Mm$ that interchanges the copy of the $\Z$-model containing $y_U$ with another copy that does not contain either $y_U$ or $F(y_U)$. There is such an automorphism by (*) in the proof of Lemma \ref{irreg-lem}. Now, this is not an automorphism of $\Pi_Ub(\Mm)$ since $F$ is one-to-one in $b(\Mm)$ (by $b(\Mm)\in\Mod(\T(S_1))$). Therefore, we assume in the rest}

Next we show that $b(\Pi_U\Mm)\ne\Pi_U b(\Mm)$ if $F(y)\notin \{\suc^k(y) : k\in\Z\}$ for some $y\in\Pi_U\omega$.  Indeed, assume the latter.  Choose an automorphism of $\Pi_U\Mm$ that interchanges the copy of the $\Z$-model containing $y$ with another copy that does not contain either $y$ or $F(y)$ and is identity on the rest. There is such an automorphism by (*) in the proof of Lemma \ref{irreg-lem}. Now, this is not an automorphism of $\Pi_Ub(\Mm)$ since $F$ is one-to-one in $b(\Mm)$ (by $b(\Mm)\in\Mod(\T(S_1))$). Therefore, we assume in the rest
\begin{description}\label{2}
	\item{(2)} $F(y)\in \{\suc^k(y) : k\in\Z\}$ for all $y\in\Pi_U\omega$.	
\end{description}
Now, (2) implies that there is a bound on ``how far $F$ can jump", i.e., there is $N_0\in\omega$ such that for all $n\in\omega$ we have
\begin{description}\label{2a}
	\item{(2a)} $F(n)=n+k$  implies $|k|<N_0$.	
\end{description}
Indeed, let $J\de\{ k\in\Z : F(n)=n+k\mbox{ for some }n\in\omega\}$ and assume that $J$ is infinite. Let $f:\omega\to J$ be a bijection, there is such a bijection because $J$ is countably infinite. For all $j\in J$ let $n_j\in\omega$ be such that $F(n_j)=n_j+j$ and let $y_i\de n_{f(i)}$ for all $i\in\omega$. Let $y \de\langle y_i : i\in\omega\rangle$. Then $F(y_U)\notin\{\suc^k(y_U) : k\in\Z\}$ because $U$ is nonprincipal. This contradicts (2), and thus $J$ is finite which implies the existence of the bound $N_0$.

Next we show that $b(\Pi_U\Mm)\ne\Pi_U b(\Mm)$ if $F$ does not agree with $\suc$ on copies of $\Z$-models in $\Pi_U\Mm$ all elements of which are in $\R$ or no elements of which are in $\R$.
Indeed, assume $\R(\suc^m(y))$ in $\Pi_U\Mm$ for all $m\in\Z$. 
There is $k\in\Z$ such that $F(y)=\suc^k(y)$, by (2). There is an automorphism $\alpha$ in $\Pi_U\Mm$ that ``shifts with $1$ step in $Y\de\{\suc^m(y) : m\in \Z\}$", i.e., $\alpha(z)=\suc(z)$ for all $z\in Y$, because $\R(z)$ for all $z\in Y$.  Now, if $F(z)\ne\suc^k(z)$ for some $z\in Y$, then $\alpha$ is not an automorphism in $\Pi_Ub(\Mm)$. So, assume that $F(z)=\suc^k(z)$ for all $z\in Y$. Now, if $k\notin\{1,-1\}$, then $Y\ne\{ F^m(y) : m\in\Z\}=\{\suc^{mk}(y) : m\in\Z\}$. However, there is an automorphism  $\beta$ of $\Pi_Ub(\Mm)$ that ``shifts  $\{ F^m(y) : m\in\Z\}$ with one step" and leaves all the other elements fixed. This $\beta$ is not an automorphism of $\Pi_U\Mm$. We show now that $F(z)=\suc^{-1}(z)$ for all $z\in Y$ cannot happen. Indeed, assume that $F(z)=\suc^{-1}(z)$ for all $z\in Y$. Then there is an ``$N_0$-long descending $F$-chain in $b(\Mm)$", i.e., there is $n\in\omega$ such that $F(k)=k-1$ for all $n-N_0\le k\le n$ in $b(\Mm)$. Then $F$ has to stay below $n$ since then on, by (2a) and $F$ being one-to-one, i.e., $F^k(n)\le n$ for all $k\in\omega$.
This again contradicts $F$ being one-to-one.
The same argument works if $\lnot\R(\suc^m(y))$ in $\Pi_U\Mm$ for all $m\in\Z$. By the above, we assume in the rest

\begin{description}\label{3}
	\item{(3)} $F(z)=\suc(z)$ for all $z\in Y\de\{\suc^k(y) : k\in\Z\}$ if $y\in\Pi_U\omega$ is such that either $\R(z)$ in $\Pi_U\Mm$ for all $z\in Y$ or $\lnot\R(z)$ in $\Pi_U\Mm$ for all $z\in Y$.
\end{description} 

\noindent
Now, (3) has implications on behavior of $F$ on long $\R$-chains or $\lnot\R$-chains in $\Mm$, as follows. Let us say that $\langle y+k : k<n\rangle = \langle y, y+1, y+2,...,y+n-1\rangle$ is an $n$-long $\R$-chain in $\Mm$ beginning with $y$ if $\R(y+k)$ in $\Mm$ for all $k<n$. The definition of a $\lnot\R$-chain is analogous. First we show the existence of a bound $N$ such that for all $\R$-chains longer than $2N$, $F$ agrees with $\suc$ on the chain, except for $N$-long chains at the beginning and at the end of the chain, and the same holds for $\lnot\R$-chains.
\begin{description}\label{3a}
	\item{(3a)} There is $N > N_0$ such that for all $\R$-chains longer than $2N$ and beginning with $y$ we have $F(\suc^k(y))=\suc^{k+1}(y)$ for all $y+N\le k\le y+n-N$ and the same holds for $\lnot\R$-chains, too.
\end{description} 
\indent Indeed, assume that there is no such bound. Then $n$ is not such a bound for any $n\in\omega$, i.e., there is an $m$-long $\R$-chain with beginning $y$ such that  $m\ge 2n$ and $F(\suc^{k}(y))\ne\suc^{k+1}(y)$ for some $y+n\le k\le y+m-n$. For each $n\in\omega$ let $y_{n}\de \suc^k(y)$ for such a chain and let $z=\langle y_{n}: n\in\omega\rangle$. Then in the ultrapower $\Pi_U\Mm$ we have $F(z)\ne\suc(z)$ while $\R(\suc^k(z))$ for all $k\in\Z$. This contradicts (3). The proof for the $\lnot\R$-chains is analogous. This completes the proof of (3a). 

\emph{From now on we assume that $N$ is as in (3a).} Next we prove that if there is an $n\ge 3N$-long $\R$-chain ending with $y-1$ and there is an $n\ge 3N$-long $\lnot\R$-chain starting with $y+1$, then the behavior of $F$ is rather close to that of $\suc$ in these chains. Namely, $F^k(o)=k$ in the interval $[y-n+N, y+n-N]$ except in $[y-N,y+N]$, and $F$ enumerates the elements of $[y-N,y+N]$. 

\begin{description}\label{3b}
	\item{(3b)} Assume that $n>3N$ and there is an $n$-long $\R$-chain in $\Mm$ ending with $y-1$ and there is an $n$-long $\lnot\R$-chain starting with $y+1$. Then $F^k(o)=k$ for all $y-n+N\le k\le y-N$ and $y+N\le k\le y+n-N$. Further, $\{ F^k(o) : y-N\le k\le y+N\}=\{ k : y-N\le k\le y+N\}$.
\end{description} 
\noindent Indeed, assume that $n$ and $y$ are as in (3b).  There is an $n\ge 2N$-long $\R$-chain beginning with $y-n$, so by (3a) we have $F(y-n+k)=y-n+k+1$ for all $y-n+N\le k\le y-N$. Let $v\de y-n+N$.  Then
\[ \tag{a} \mbox{$F(v+k)=v+k+1$ for all $k\le n-2N$.}\]
 Then $F(w)\notin\{ k : v< k \le v+n-2N\}$ for all $w<v$ since $F$ is one-to-one by $b(\Mm)\models\T(S_1)$. By $n-2N\ge N>N_0$ and (2a) then $F(w) \le v$ for all $w<v$ and hence $F$ enumerates $[0,v]$, i.e., 
 \[ \tag{b} \mbox{$\{ F^k(o) : k\le v\}=\{ k : k\le v\}$.}\]
  There is $m\in\omega$ such that $v=F^m(o)$, by (1). As before, by (2a) and (a) we have that $m\le v$  and then $m=v$  by (b). Thus, $F^{v}(o)=v$ and by (a) we have $F^k(o)=k$ for all $v\le k\le y-N$. The rest of (3a) can be obtained similarly.
\bigskip

We are ready to show $b(\Pi_U\Mm)\ne\Pi_Ub(\Mm)$, finishing the proof of Theorem \ref{postest-cor}. Let $X$ be the infinite set where $S_0$ and $S_1$ differ. Then $X$ is disjoint from $S_0$ and $X\subseteq S_1$, by definition.  Let  $x_n$ denote the $n$th member of $X$ according the natural ordering of $\omega$.  Then $\lnot\R(x_n)$ in $\Mm$ by $x_n\notin S_0$ and the definition of $\R$ in $\Mm$. Also, $\R(x_n-k-1)$ and $\lnot\R(x_n+k)$ for all $k<n$, because the $0,1$-sequences between two $x$s are laid by alphabetical order, thus before the $n$th  $x\in X$ there are $n$ many $1$s and after it there are $n+1$ many $0$s. Let $x\de\langle x_n : n\in\omega\rangle$. Then $x_U$ is contained in $\Pi_U\Mm$ in a copy of the $\Z$-model $\langle \Z,\{k : k<0\}\rangle$, i.e., all members of the $\Z$-model below $x_U$ are in $\R$, and no member after $x_U$, including $x_U$ is in $\R$. 

How does the set $C\de\{\suc^k(x_U) : k\in\Z\}$ look like in $\Pi_Ub(\Mm)$? Note that we cannot assume $F=\suc$ and $o=0$ in $b(\Mm)$. Thus,  for example, we cannot infer $\R(x_n)$ in $b(\Mm)$ from $x_n\in S_1$. However, we can use our assumptions (1)-(3) and their implications. Especially, we can use (3b).
Let $n\ge 3N$,  where $N$ is the bound in (3b). We have seen in the previous paragraph that, in $\Mm$, the assumptions hold for $y=x_n$. By (1), the definition of $S_1$, and (3b) then  $\R(F^k(o))$ for $x_n-n+N\le k\le x_n-N$ and $\lnot\R(F^k(o))$ for $x_n+N\le k\le x_n+n-N$, in $b(\Mm)$. Also, by (3b) we get that $F$ agrees with $\suc$ ``below" $\suc^{-N}(x_U)$ and ``above" $\suc^{N}(x_U)$, in $\Pi_Ub(\Mm)$. Further, $F$ enumerates the interval $I\de[\suc^{-N}(x_U),\suc^N(x_U)]$. However, there is a difference between $\Pi_U\Mm$ and $\Pi_Ub(\Mm)$ concerning $I$. Namely, in $\Pi_U\Mm$ exactly  $N$ elements of $I$ are in $\R$ because $\lnot\R(x_n)$ in $\Mm$. At the same time, due to the definition of $S_1$, by (1) we get $\R(F^w(o))$ for all $w\in X$. Hence,  exactly $N+1$ elements are in $\R$ in the corresponding intervals in $b(\Mm)$, so exactly $N+1$ elements of $I$ are in $\R$, in $\Pi_Ub(\Mm)$. 

For all $n\in\omega$ let $y_n\in\omega$ be similar to $x_n$ in that $\lnot\R(y_n)$, there is an $n$-long $\R$-chain ending with $y_n-1$, there is an $n$-long $\lnot\R$-chain starting with $y_n+1$, and such that neither $y_n$ nor any element of these chains belong to $X$. There are such $y_n$s by the construction of $S_0, S_1$. Let $y\de\langle y_n : n\in\omega\rangle$. Then there is an automorphism in $\Pi_U\Mm$ that interchanges $x_U$ with $y_U$.  
We will show that there is no automorphism in $\Pi_Ub(\Mm)$ that interchanges $\suc^{-N}(x_U)$ and $\suc^{-N}(y_U)$. Indeed, such an automorphism has to be a bijection between the intervals $I$ and $J$ because it can be seen that $F$ enumerates $J\de[\suc^{-N}(y_U),\suc^N(y_U)]$ in $\Pi_Ub(\Mm)$ and $F$ agrees with $\suc$ outside $J$. We have seen that there are $N+1$ elements of $I$ that are in $\R$ in $\Pi_Ub(\Mm)$. It can be seen just the same way that there are only $N$ elements of $J$ because $\lnot\R(y_U)$ in $\Pi_Ub(\Mm)$. Therefore, no bijection between $I$ and $J$ can preserve $\R$. The proof of Theorem \ref{postest-cor} is complete. 
\qed\bigskip

We close the paper with some implications of the results for questions raised in the wider literature.\bigskip

Glymour~\cite{Gly} raises an interesting question about definitional equivalence. The common understanding is that definitionally equivalent theories have essentially the same content and we would think that all important properties are shared by them. 
Theorem~\ref{main-thm} implies that a property of a theory is preserved by definitional equivalence when it can be expressed in terms of universes, isomorphisms and ultraproducts of its models. Therefore, having a one-element model, having only finite models, being categorical in a power or being complete are preserved by classical definitional equivalence (since two models are elementarily equivalent if and only if they have isomorphic ultrapowers). Glymour \cite[p.296]{Gly} conjectures that also the model class being closed under substructures, the model class being closed under unions of chains, and having an equational axiomatization are preserved. We now show that neither one of these three properties is preserved by definitional equivalence.

Indeed, let $\T_1$ be the empty theory on the language with one constant symbol $c$. Let $\T_2$ be the definitional extension of $\T_1$ with
$\forall x(R(x)\leftrightarrow [\exists yz(y\ne z)\land x=c]$. Then $\Mod(\T_1)$ is closed under taking substructures but $\Mod(\T_2)$ is not. The counterexample to preservation of unions of chains is similar in spirit. Let $\T_1$ be the empty theory on the language whith a binary relation symbol $\le$. Let $\T_2$ be the definitional extension of $\T_1$ with defining $R$ to be the set of $\le$-minimal elements when there is a $\le$-maximal element and $R$ is the empty set when there is no $\le$-maximal element (i.e., $\forall x[R(x)\leftrightarrow(\exists y\forall z(z\le y)\land \forall z(z\le x))]$). Clearly, $\T_1$ is closed under taking unions of chains. However, $\T_2$ is not closed under taking unions of chains, as the following models show. For each natural number $n$ let $\Mm_n$ have the set of natural numbers smaller than $n$ as universe, let $\le$ be the ``smaller" relation and let only $0$ be in relation $R$. Then each $\Mm_n$ is a model of $\T_2$ but their union is not a model of $\T_2$ since it does not have a maximal element yet $R$ is nonempty in it. For showing that having an equational axiomatization is not preserved by definitional equivalence, 
one could take groups as counterexamples, this is mentioned in \cite[p.56]{HMT71}. Indeed, let  $\T_1$ be the class of semigroups in which inverses exist and let $\T_2$ be its extension with the inverse operation and the zero element as constant. Then $\T_1$ does not have a universal axiomatization because its model class is not closed under subalgebras, while $\T_2$ is an equational class.

It is known that definitionally equivalent theories have isomorphic Linden\-baum--Tarski formula-algebras, they only differ from each other in what definable properties they take to be as basic ones. The proofs above show that this latter choice can influence the existence of axiom systems of given forms.
For example, being substructure is not preserved by definitional expansion because in this notion the basic relations are treated differently from the rest, namely being a substructure is formulated in terms of basic relations only. Similarly for homomorphism, union, etc. However, being an elementary substructure is preserved by definitional expansion because in the definition of elementary substructure all definable relations are treated alike (and indeed, this notion can be characterized by means of isomorphisms and ultraproducts as follows: $\Nn$ is an elementary substructure of $\Mm$ if and only if $N\subseteq M$ and there is an ultrafilter $U$ such that $\Pi_U\Nn$ is isomorphic to $\Pi_U\Mm$ via an isomorphism that is identity on the diagonal image of $N$ in $\Pi_U\Nn$).  
\bigskip

The following corollary of Theorem~\ref{main-thm} states that an associated structure to be
defined below, namely the concrete ultracategory of a theory, is an
invariant characteristic to definitional equivalence of first-order
logic theories.

By a \emph{concrete ultracategory}, we mean a triple  $(C,F,p)$ where 
$(C,F)$ is a concrete category\footnote{For the notions of a concrete category and a concrete
functor see \cite[Chap.5]{AHS}.},
and the additional structure $p$ is a system of infinitary functions $\langle p_U : U \mbox{ an
ultrafilter}\rangle$ on $Ob(C)$ such that if $U$ is an ultrafilter on the set $I$
then $F(p_U(m_i)_{i\in I})=\Pi_UF(m_i)$ for all $m:I\to Ob(C)$.
A functor between two ultracategories $(C,F,p)$ and $(C',F',p')$ is a concrete functor between $(C,F)$ and $(C',F')$
that preserves all the functions $p_U$.  Two concrete ultracategories
are isomorphic if there is a functor between them that is a category
theoretical isomorphism. 

Let $\T$ be a theory. Its \emph{concrete ultracategory} is $(C,F,p)$
where $(C,F)$ is $\Cat^{iso}(\T)$ with the natural forgetful
functor, and for all ultrafilters $U$ on $I$ and all systems $(\Mm_i)_{i\in I}$ we have
$p_U((\Mm_i)_{i\in I})=\Pi_U\Mm_i$. Notice that an isomorphism between the ultracategories of two theories preserves only the universes of the models (through the forgetful functors) and the behaviour of isomorphisms and ultraproducts as functions on $\Cat^{iso}(\T)$.

\begin{thm}\label{ultracat-thm} Two first-order logic theories are
definitionally equivalent if and only if their concrete
ultracategories are isomorphic.
\end{thm}

\noindent{\bf Proof.} 
This is just a reformulation of
Theorem~\ref{main-thm}. 
\qed\bigskip

We note that one can define the concrete ultracategory of a theory
to contain all elementary embeddings in place of all isomorphisms
only, as is usual.
Theorem~\ref{ultracat-thm} is true with this modified definition,
too. The reason is that elementary embeddings are preserved by
definitional equivalence.
\bigskip

\noindent{\bf Remark 7.}(connection with Stone duality)
Halvorson~\cite[section 7]{Halwhat} proposes the programme to investigate what structure a model class naturally 
has. This program involves to endow the model class of a theory in such a way that from this structure on the model class, the theory can be recovered up to definitional equivalence. Theorem~\ref{ultracat-thm} above offers an answer, namely concrete ultracategory of a theory. In category theoretical logic, Makkai \cite[Theorem 4.1]{Makkai87} offers the notion of (abstract) ultracategory and Awodey and Forssell \cite{Awo} offer the notion of topological groupoid in place of our concrete ultracategory. These three structures are quite similar to each other, so there seems to be a convergence here in finding a natural structure on the model classes. Unlike our concrete ultracategory, Makkai's ultracategory and Awodey and Forssell's topological groupoids characterize first-order theories only up to many-sorted definitional equivalence, which is weaker than classical definitional equivalence. 
Halvorson~\cite{Halwhat} points out the connection of his programme with generalizing Stone duality from propositional logic to predicate logic. We believe that a full-fledged Stone duality can be based on Theorem~\ref{ultracat-thm} above. See also \cite{Makkai85, Makkai87, harnik} and
\cite[p.576]{BH16}.
\comment{The above theorem is analogous to Michael Makkai's ultracategory
theorem which is available at present, as far as we know, only in
category theoretical form \cite[Theorem 4.1]{Makkai87}.
Makkai's theorem is about the more general and more refined
definitional equivalence which allows one to define new universes
(and thus automorphism groups are preserved only up to isomorphism),
while Theorem~\ref{ultracat-thm} is about the simpler classical
definability. See \cite{Makkai85, Makkai87, harnik} and
\cite[p.576]{BH16}.} \qquad $\Box$
\bigskip

Definability theory is used quite extensively in recent philosophy of science papers to investigate what symmetries tell about theories and how to compare ``structure", see, for example, \cite{Bsym, BMW, Halvorson, HudetzDCE}. When one theory is an expansion of the other, there is a natural functor between their model categories. This is the ``reduct-formation" functor denoted by $\Pi$ in \cite[above Example 9]{Bsym}. It is shown in \cite{Bsym} that the question investigated in the present paper gets rather nice answers in this special case. We now show how one of the attractive theorems in \cite{Bsym} follows from Theorem \ref{main-thm}. In fact, Theorem \ref{main-thm} in the present paper is a generalization of \cite[Corollary 2]{Bsym} to the general case concerning two arbitrary theories.

\begin{cor}\label{cor1}(Corollary 2 in \cite{Bsym})
Let $\T^{+}$ be an expansion of $\T$. Then $\T^{+}$ is definitionally equivalent to $\T$ if and only if the reduct-formation functor $\Pi$ is an equivalence between their model iso-categories.
\end{cor}

\noindent{\bf Proof.} 
The reduct-formation functor $\Pi$ is a concrete functor and it always preserves isomorphisms and ultraproducts ``forwards", i.e., from $\T^{+}$ to $\T$. It is a bijection up to isomorphism if and only if it is a bijection because the range of $\Pi$ is always closed under isomorphisms. Thus if $\Pi$ is a category theoretical equivalence then each model of $\T$ has a unique expansion in $\Mod(\T^{+})$, therefore $\Pi$ preserves isomorphisms and ultraproducts also backwards. Thus if $\Pi$ is a category theoretical equivalence then it satisfies (i) and (ii) in Theorem \ref{main-thm}, hence $\T$ and $\T^+$ are definitionally equivalent. The other direction is easy. 
	\qed\bigskip

Categorical equivalence of theories is investigated in \cite{BH16} as a weaker form of definitional equivalence. Two theories are defined to be \emph{categorically equivalent} iff there is a categorical equivalence between their model categories. It is shown in \cite{BH16} that categorical equivalence, many-dimensional (Morita) equivalence and definitional equivalence are strictly stronger in this order. The question naturally arises about how ``large" the gaps between them are, under what additional properties these are the same. 

According to Corollary \ref{cor1}, the reduct-formation functor $\Pi$ bridges the gap between definitional equivalence and categorical equivalence between a theory and its expansion. It is asked in \cite[below Corollary 2]{Bsym} what special property $\Pp$ of $\Pi$ allows it to fill the gap between categorical and definitional equivalence of theories. Theorem \ref{main-thm} gives an answer to this question. The answer it offers is that this special property $\Pp$ of $\Pi$ is that it is a concrete functor which preserves ultraproducts in both directions when it is an equivalence.

Question 2 in \cite{BDis17} asks for an additional property $\Pp$ of functors
such that two theories are definitionally equivalent iff there is a category
theoretical equivalence between their model categories which has
property $\Pp$. This question is also mentioned in \cite[Note
23]{Weatherall}, where it is written: ``It is not known how much
weaker categorical isomorphism is than definitional equivalence, or
Morita equivalence, which is a weakening of definitional equivalence
that allows one to define new sorts." Now, Corollary \ref{dist-cor} below says, roughly, that 
categorical equivalence is just as much weaker than
definitional equivalence as it misses how ultraproducts behave and what the universes of models as well as the set theoretical contents of morphisms are.
In other words, two theories are definitionally equivalent if and only if there is an equivalence between their model categories which is a concrete isomorphism  and preserves ultraproducts. We note that \cite[Theorem 3]{HudetzDCE}
gives an answer to the above questions that is different
in spirit from our Corollary~\ref{dist-cor}.

\begin{cor}\label{dist-cor}
Two theories  $\T_1$ and $\T_2$ are definitionally equivalent if and
only if there is a  concrete ultraproduct-preserving functor $F$ that is an
equivalence between $\Cat(\T_1)$ and $\Cat(\T_2)$.
\end{cor}

\bigskip
Ultraproducts are intimately connected to first-order logic. It
would be interesting to see whether analogous theorems hold for
other languages where ultraproducts can be omitted or replaced with
some other additional structure. Laurenz Hudetz \cite{HudetzDCE,
Hudetz} contain interesting generalizations and results in the
direction of broadening definability theory in order to be more
applicable in philosophy of science. These results may be used
perhaps to get an analogue of Theorem~\ref{main-thm} in which
ultraproducts do not occur.

\appendix
\section{Appendix}

The following generalized version of Beth's theorem is well-known as
folklore. Both \cite{vBP} and \cite{Makkai87} use this generalized version of Beth's theorem without proof.  Since Theorem \ref{main-thm} relies heavily on this folklore theorem, here we give a short proof for it. For simplicity, we assume that we have only relation symbols.

\begin{thm}\label{genBeth-thm} Assume that $\T$ is a theory on the language $\Sigma\cup\RR$ and the $\Sigma$-reduct of each model of $\T$ has at most one extension to a model of $\T$. Then each element of $\RR$ is explicitly definable in $\T$ by a $\Sigma$-formula.
\end{thm}

\noindent{\bf Proof.} 
Let $\T'$ denote the theory $\T$ where each relation symbol $R\in\RR$ is replaced by a new relation symbol $R'$ not occurring in the language of $\T$ (and having the same arity). Then $\T\cup\T'\models \forall\bar{x}[R(\bar{x})\leftrightarrow R'(\bar{x})]$ for all $\R\in\RR$, since the $\RR$-free reduct of each model of $\T$ has at most one expansion to a model of $\T$. Let $R\in\RR$ be arbitrary. By the compactness theorem, there is a finite subset $\T_0$ of $\T$ such that $\T_0\cup\T_0'\models \forall\bar{x}[R(\bar{x})\leftrightarrow R'(\bar{x})]$.  Therefore, $R$ has to occur in $\T_0$, since otherwise both the empty set and the biggest relation of the same rank as $R$ can be chosen in a model to satisfy $\T_0$.  Since $\T_0$ is finite, it contains only finitely many elements from $\RR$, let the set of these elements be $\RR_0\de\{ R_1,\dots,R_n\}$, and we may assume $R_1$ is $R$. By the usual Beth's theorem, there is a formula $\varphi_R$ on the language $\Sigma\cup\{ R_2,\dots, R_n\}$ which defines $R$ in $\T_0$.  Now, let $\T_1$ be the theory we obtain from $\T_0$ by replacing $R$ in it everywhere with $\varphi_R$. Then $\T_1$ follows from $\T_0$, only $R_2,\dots, R_n$ occur in $\T_1$ and $\T_1\cup\T_1'\models \forall\bar{x}[R_2(\bar{x})\leftrightarrow R_2'(\bar{x})]$. By the usual Beth's theorem, there is a formula $\varphi_{R1}$ on the language $\Sigma\cup\{ R_3,\dots, R_n\}$ which defines $R_1$ in $\T_1$. And so on. At the end we get $\T_{n-1}$ on the language $\Sigma\cup\{ R_n\}$ and a formula $\varphi_{Rn}$ on the language $\Sigma$ which defines $R_n$ in $\T_{n-1}$. Let $\psi_n$ be $\varphi_{Rn}$, let $\psi_{n-1}$ be the formula we get from $\varphi_{Rn-1}$ by replacing $R_n$ in it by $\psi_n$, etc. Then $\psi_1$ is in the language $\Sigma$ which defines $R$ in $\T_0\subseteq \T$. 
\qed\bigskip

\noindent{\bf Acknowledgements} We are indepted to the two anonymous referees for their very useful feedbacks.

\vspace*{10pt}

\comment{\address{ALFR\'ED R\'ENYI INSTITUTE OF MATHEMATICS\\
BUDAPEST, RE\'ALTANODA st.\ 13-15, H-1053 HUNGARY\\
{\it E-mail}: {andreka.hajnal, madarasz.judit, nemeti.istvan, szekely.gergely}@renyi.hu\\
UNIVERSITY OF PUBLIC SERVICE\\
BUDAPEST, 2 LUDOVIKA square, H-1083 HUNGARY}
\clearpage}
\noindent ALFR\'ED R\'ENYI INSTITUTE OF MATHEMATICS\\
		BUDAPEST, RE\'ALTANODA st.\ 13-15, H-1053 HUNGARY\\
		{\it E-mail}: {andreka.hajnal, madarasz.judit, nemeti.istvan, szekely.gergely}@renyi.hu\\
		UNIVERSITY OF PUBLIC SERVICE\\
		BUDAPEST, 2 LUDOVIKA square, H-1083 HUNGARY

\end{document}